\documentclass[11pt,twoside]{article}

\usepackage{a4wide}
\usepackage{amsfonts}
\usepackage{amssymb}
\usepackage{amsmath}
\usepackage{graphicx}

\newcommand{\ignore}[1]{}

\def\@begintheorem#1#2{\par\bgroup{\sc #1\ #2. }\it\ignorespaces}
\def\@opargbegintheorem#1#2#3{\par\bgroup{\sc #1\ #2\ (#3). } \it\ignorespaces}
\def\@endtheorem{\egroup}
\newtheorem{theorem}{Theorem}[section]
\newtheorem{corollary}[theorem]{Corollary}
\newtheorem{lemma}[theorem]{Lemma}

\newtheorem{example}[theorem]{Example}
\newtheorem{proposition}[theorem]{Proposition}
\newtheorem{definition}[theorem]{Definition}
\newcommand{\bt}[1]{\begin{theorem}\label{#1}}
\newcommand{\bc}[1]{\begin{corollary}\label{#1}}
\newcommand{\bl}[1]{\begin{lemma}\label{#1}}
\newcommand{\be}[1]{\begin{example}\label{#1}}
\newcommand{\bp}[1]{\begin{proposition}\label{#1}}
\newcommand{\ba}[1]{\begin{algorithm}\rm\label{#1}}
\newcommand{\bd}[1]{\begin{definition}\rm\label{#1}}{\normalsize }
\newcommand{\bpr}{\noindent {\em Proof. }}
\newcommand{\et}{\end{theorem}}
\newcommand{\ec}{\end{corollary}}
\newcommand{\el}{\end{lemma}}
\newcommand{\ee}{\end{example}}
\newcommand{\ep}{\end{proposition}}
\newcommand{\ed}{\end{definition}}
\newcommand{\epr}{{\ \vbox{\hrule\hbox{%
\vrule height1.3ex\hskip0.8ex\vrule}\hrule}}\\\par}

\def\R{\mathbb{R}}
\def\Z{\mathbb{Z}}

\begin{document}

\title{\bf Optimization over Degree Sequences of Graphs}

\author{
Gabriel Deza
\thanks{\small University of Toronto. Email: gabe.deza@mail.utoronto.ca}
\and
Shmuel Onn
\thanks{\small Technion - Israel Institute of Technology.
Email: onn@ie.technion.ac.il}
}
\date{}

\maketitle

\begin{abstract}
We consider the problem of finding a subgraph of a given graph minimizing the sum of
given functions at vertices evaluated at their subgraph degrees. While the problem
is NP-hard already for bipartite graphs when the functions are convex on one side
and concave on the other, we show that when all functions are convex, the problem
can be solved in polynomial time for any graph. We also provide polynomial time solutions
for bipartite graphs with one side fixed for arbitrary functions, and for arbitrary graphs
when all but a fixed number of functions are either nondecreasing or nonincreasing. We
note that the general factor problem and the (l,u)-factor problem over a graph are special
cases of our problem, as well as the intriguing exact matching problem. The complexity
of the problem remains widely open, particularly for arbitrary functions over complete graphs.

\end{abstract}

\section{Introduction}

The {\em degree sequence} of a simple graph $G=(V,E)$ with $V=[n]:=\{1,\dots,n\}$
is the vector $d(G)=(d_1(G),\dots,d_n(G))$, where $d_i(G):=|\{e\in E:i\in e\}|$
is the degree of vertex $i$ for all $i\in[n]$. Degree sequences have been studied
by many authors, starting from their celebrated effective characterization by
Erd\H{o}s and Gallai \cite{EG}, see e.g. \cite{EKM} and the references therein.

\vskip.2cm
In this article we are interested in the following discrete optimization problem.

\vskip.2cm\noindent{\bf Optimization over Degree Sequences.}
Given a graph $H$ on $[n]$, and for $i=1,\dots,n$ functions
$f_i:\{0,\dots,d_i(H)\}\rightarrow\Z$, find a subgraph $G\subseteq H$
on $[n]$ minimizing $\sum_{i=1}^n f_i(d_i(G))$.

\vskip.2cm
A special case of this is the {\em general factor problem} \cite{Cor},
which is to decide, given a graph $H$ on $[n]$ and subsets $B_i\subseteq\{0,1,\dots,d_i(H)\}$
for $i\in[n]$, if there is a $G\subseteq H$ with $d_i(G)\in B_i$ for all $i$,
and find one if yes. Indeed, for each $i$ define a
function $f_i$ by $f_i(x):=0$ if $x\in B_i$ and $f_i(x):=1$ if $x\notin B_i$.
Then the optimal value of the degree sequence problem is zero if and only if there is a factor, in which case any optimal graph
$G$ is one. An even more special case is the well studied {\em $(l,u)$-factor problem}
introduced by Lov\'asz \cite{Lov}, where each $B_i=\{l_i,\dots,u_i\}$ is an interval.
This reduces to our degree sequence problem even with convex functions, with
$f_i(x):=l_i-x$ if $0\leq x\leq l_i$, $f_i(x):=0$ if $l_i\leq x\leq u_i$,
and $f_i(x):=x-u_i$ if $u_i\leq x\leq d_i(H)$.

\vskip.2cm
Another interesting special case is when $H=K_n$, that is, when the optimization is over any graph
$G$ on $[n]$. This special case was recently shown in \cite{DLMO} to be solvable in polynomial time
when all the functions are the same, $f_1=\cdots=f_n$, extending an earlier result of \cite{PPS}.

\vskip.2cm
We note right away that the degree sequence problem in its full generality is hard.
\bp{hard}
Deciding if the optimal value in our problem is zero is NP-complete already:
\begin{enumerate}
\item
when $f_1=\cdots=f_n=f$ are identical, with $f(0)=f(3)=0$ and $f(i)=1$ for $i\neq0,3$.
\item
when $H=(I,J,E)$ is bipartite and $f_i$ is convex for all $i\in I$ and concave for all $i\in J$.
\end{enumerate}
\ep

In contrast, for convex functions, we prove the following statement:
\bt{convex}
The optimization problem over degree sequences can be solved in polynomial time for any
given $n$, any given graph $H$ on $[n]$, and any given convex functions $f_1,\dots,f_n$.
\et
This in particular allows to solve the $(l,u)$-factor problem via the reduction described above.

\vskip.2cm
Next, for unbalanced bipartite graphs, that is, with one side fixed and small, we prove:
\bt{unbalanced_bipartite}
For any fixed $r$, the optimization problem over degree sequences can be solved in polynomial
time for any given bipartite graph $H\subseteq K_{r,n-r}$ and any given functions $f_1,\dots,f_n$.
\et

Finally, when most functions are all nondecreasing or all nonincreasing we prove:
\bt{fixed_monotone}
For any fixed $r$, our problem is polynomial time solvable for any graph $H$ and any functions
$f_1,\dots,f_r$ when $f_{r+1},\dots,f_n$ are either all nondecreasing or all nonincreasing.
\et

\vskip.2cm
Optimization over degree sequences naturally extends to hypergraphs, see e.g. \cite{Bil,DLMO}.
But this is harder, since deciding if the optimal value given convex functions $f_i(x)=(x-d_i)^2$
is zero, is equivalent to deciding if there is a hypergraph with degrees $d_1,\dots,d_n$, which
was recently shown in \cite{DLMO} to be NP-complete, solving a long standing open problem from \cite{CKS}.

The general factor problem can be solved in polynomial time when no $B_i$ has a gap longer than
$1$ (for instance, $\{0,3,4,6\}$ has longest gap $2$), see \cite{Cor}. This implies that also the
degree sequence problem can be solved for the functions associated to such $B_i$ as described above.

However, beyond the new results and surveyed results above, the degree sequence problem is
still wide open. It would be interesting to determine the complexity of the  problem for
various classes of input graphs $H$ and various classes of input functions $f_i$. Particularly
intriguing is the special case of $H=K_n$, when the optimization is over any graph $G$ on $[n]$,
where it may be that the problem can be solved in polynomial time for any given functions.

\vskip.2cm
The article is organized as follows. In Section \ref{hem} we prove Proposition \ref{hard}
and also show a simple reduction of the well known and intriguing {\em exact matching problem}
to our problem. In Section \ref{cf} we consider convex functions and prove Theorem \ref{convex}.
In Section \ref{ubg} we consider unbalanced bipartite graphs, that is, with one side fixed
and small, and prove Theorem \ref{unbalanced_bipartite}. In Section \ref{mf} we consider
nondecreasing or nonincreasing functions and prove Theorem \ref{fixed_monotone}.

\section{Hardness and exact matching}
\label{hem}

We begin by showing that the optimization problem over degree sequence is generally hard.

\vskip.2cm\noindent{\bf Proposition \ref{hard}}
{\it
Deciding if the optimal value in our problem is zero is NP-complete already:
\begin{enumerate}
\item
when $f_1=\cdots=f_n=f$ are identical, with $f(0)=f(3)=0$ and $f(i)=1$ for $i\neq0,3$.
\item
when $H=(I,J,E)$ is bipartite and $f_i$ is convex for all $i\in I$ and concave for all $i\in J$.
\end{enumerate}
}
\bpr
\begin{enumerate}
\item
The NP-complete {\em cubic subgraph problem} \cite{GJ} is the case of the factor problem with
$B_i=\{0,3\}$ for all $i$. Defining $f_i$ as described above, our problem is NP-complete too.
\item
The general factor problem is NP-complete already when $H=(I,J,E)$ is bipartite
with maximum degree $3$ and $B_i$ equals $\{1\}$ for all $i\in I$
and $\{0,3\}$ for all $i\in J$, see \cite{Cor}. Define
$$
f_i(x)\ :=\ (x-1)^2
\quad i\in I,
\quad\quad\quad\quad
f_i(x)\ :=\ x(3-x)
\quad i\in J\quad.
$$
Then $f_i$ is convex for all $i\in I$ and concave for all $i\in J$, and the optimal value of
the degree sequence problem is zero if and only if there is a factor, proving the claim.
\end{enumerate}
\epr

The {\em exact matching problem} is the following. Fix a positive integer $r$. Given
$b\in\Z_+^r$ and a labeling $c:[n]\times[n]\rightarrow[r]$ of the edges of $K_{n,n}$,
decide if $K_{n,n}$ has a perfect matching $M\subset[n]\times[n]$ such that
$|\{e\in M:c(e)=k\}|=b_k$ for all $k$, and find one if yes. This problem has randomized
algorithms \cite{MVV,Onn} but its deterministic complexity is open already for $r=3$.

We note that it is a simple special case of the general factor problem and hence of the degree
sequence problem over subgraphs. Let $b$, $c$ be given. Define a graph $H=(V,E)$ by
$$V\ :=\ \{u_1,\dots,u_n\}\ \uplus\ \{v_1,\dots,v_n\}
\ \uplus\ \{w_{i,j}\, :\, 1\leq i,j\leq n\}\ \uplus\ \{x_1,\dots,x_r\}\ ,$$
$$E\ :=\ \biguplus_{1\leq i,j\leq n}
\left\{\{u_i,w_{i,j}\},\{w_{i,j},v_j\},\{w_{i,j},x_{c(i,j)}\}\right\}\ .$$
Define sets
$$B_{u_i}:=B_{v_j}:=\{1\},\ \ B_{w_{i,j}}:=\{0,3\},
\ \ B_{x_k}:=\{b_k\},\quad 1\leq i,j\leq n,\ \ 1\leq k\leq r$$
and corresponding convex and concave functions
$$f_{u_i}:=f_{v_j}:=(x-1)^2,\ \ f_{w_{i,j}}:=x(3-x),
\ \ f_{x_k}:=(x-b_k)^2,\quad 1\leq i,j\leq n,\ \ 1\leq k\leq r\ .$$
Then it is not hard to verify that there exists an exact matching
if and only if the optimal objective value of the degree sequence problem
over $H$ and the $f_i$ is zero, in which case an optimal graph $G\subseteq H$
is also a factor and $M:=\{(i,j):d_{w_{i,j}}(G)=3\}$ is an exact matching.

\section{Convex functions}
\label{cf}

Here we consider the situation where all the given functions are convex.
\vskip.2cm\noindent{\bf Theorem \ref{convex}}
{\it
The optimization problem over degree sequences can be solved in polynomial time for any
given $n$, any given graph $H$ on $[n]$, and any given convex functions $f_1,\dots,f_n$.
}

\vskip.2cm
\bpr
Let the given graph be $H=([n],E)$ and let $d(H)$ be its degree sequence.
Construct a graph $L=(V,F)$ with $8|E|$ vertices and $4|E|+2\sum_{i=1}^n d_i^2$ edges
as follows. For $i\in[n]$ let $\delta(i)=\{e\in E : i\in e\}$.
Introduce the $4d_i$ vertices $u_i^e,v_i^e$ for $e\in\delta(i)$ and $x_i^k,y_i^k$ for $k\in[d_i]$.
Introduce the $2d_i^2$ edges $\{u_i^e,x_i^k\}$ ,$\{v_i^e,y_i^k\}$ for $k\in[d_i]$ and
$e\in\delta(i)$. Introduce the $d_i$ edges $\{x_i^k,y_i^k\}$ for $k\in[d_i]$.
For $e=\{i,j\}\in E$ introduce the two edges $\{u_i^e,u_j^e\}$ and $\{v_i^e,v_j^e\}$.

Next define a cost function $c:F\rightarrow\Z$ on the edges of $L$ as follows.
For $i\in[n]$ and $k\in[d_i]$ let $c(\{x_i^k,y_i^k\}):=f_i(k)-f_i(k-1)$.
Let the costs of all other edges in $F$ be zero.

We now show that the degree sequence problem over $H$ reduces to the minimum cost
perfect matching problem over $L$ and therefore can be solved in polynomial time.

Denote for brevity $c_i^k:=c(\{x_i^k,y_i^k\})$. Then for each $i\in[n]$ we have
$c_i^1\leq\cdots\leq c_i^{d_i}$ since for each $1\leq k< d_i$ the convexity
of $f_i$ implies $2f_i(k)\leq f_i(k-1)+f_i(k+1)$ and therefore
$$c_i^k\ =\ f_i(k)-f_i(k-1)\ \leq\ f_i(k+1)-f_i(k)\ =\ c_i^{k+1}\ .$$

Let $f^*$, $c^*$ be the optimal values of the degree sequence
and matching problems respectively.

\begin{figure}[hbt]
\centerline{\includegraphics[scale=0.8]{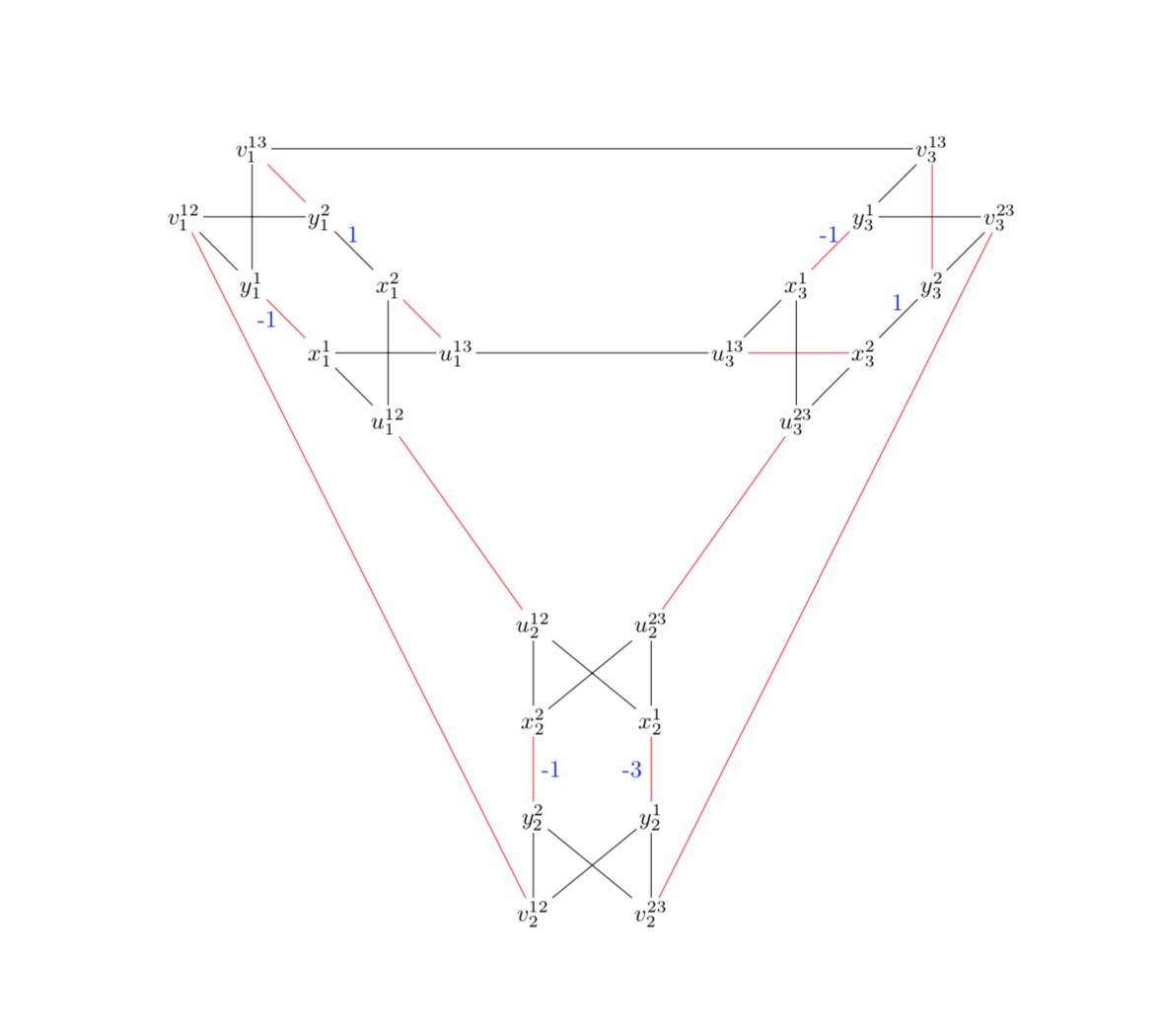}}
\caption{The graph $L$ of Example \ref{example_matching}
with blue costs and red minimum perfect cost matching}
\label{figure1}
\end{figure}

Consider first a $G\subseteq H$ optimal for the degree sequence problem with $G=([n],N)$.
Define a perfect matching $M\subset F$ of $L$ as follows. For each $e=\{i,j\}\in N$
include in $M$ the two edges $\{u_i^e,u_j^e\}$ and $\{v_i^e,v_j^e\}$. For each $i\in[n]$
include in $M$ the edges $\{x_i^k,y_i^k\}$ for $k=1,\dots,d_i(G)$. Next, match the
$d_i(H)-d_i(G)$ unmatched vertices among the $u_i^e$ for $e\in\delta(i)$ with the
unmatched vertices $x_i^k$ for $d_i(G)<k\leq d_i(H)$, and match the
$d_i(H)-d_i(G)$ unmatched vertices among the $v_i^e$ for $e\in\delta(i)$ with the
unmatched vertices $y_i^k$ for $d_i(G)<k\leq d_i(H)$.
Now we have the inequality
\begin{eqnarray}\label{1}
c^*\ \leq\ \sum_{i=1}^n\sum_{k=1}^{d_i(G)}c_i^k
& = & \sum_{i=1}^n\sum_{k=1}^{d_i(G)}\left(f_i(k)-f_i(k-1)\right)\\
\nonumber
& = & \sum_{i=1}^n\left(f_i(d_i(G))-f_i(0)\right)\ =\ f^*-\sum_{i=1}^n f_i(0)\ .
\end{eqnarray}
Next, consider a minimum cost perfect matching $M\subset F$ and define a subgraph
$G\subseteq H$ with $G=([n],N)$ where $N:=\{e=\{i,j\}\in E\, :\, \{u_i^e,u_j^e\}\in M\}$.
Then $d_i(G)$ vertices among the $u_i^e$ for $e\in\delta(i)$ are matched in $M$ with other
vertices of the form $u_j^h$ and therefore the other $d_i(H)-d_i(G)$ vertices among $u_i^e$
for $e\in\delta(i)$ are matched with vertices of the form $x_i^k$. So for some subset
$K_i\subseteq[d_i(H)]$ with $|K_i|=d_i(G)$ the vertices $x_i^k$ for $k\in K_i$
are matched with the corresponding vertices $y_i^k$. Since $c_i^1\leq\cdots\leq c_i^{d_i}$,
we now have the inequality
\begin{eqnarray}\label{2}
c^*\ =\ \sum_{i=1}^n\sum_{k\in K_i}c_i^k
\ \geq\ \sum_{i=1}^n\sum_{k=1}^{d_i(G)}c_i^k
& = & \sum_{i=1}^n\sum_{k=1}^{d_i(G)}\left(f_i(k)-f_i(k-1)\right)\\
\nonumber
& = & \sum_{i=1}^n\left(f_i(d_i(G))-f_i(0)\right)\ \geq\ f^*-\sum_{i=1}^n f_i(0)\ .
\end{eqnarray}
Equations \eqref{1} and \eqref{2} together imply that $c^*=f^*-\sum_{j=1}^n f_j(0)$
and so equality holds all along in both equations. This implies in particular
that the subgraph $G\subseteq H$ obtained from a minimum cost
perfect matching as above is optimal for the degree sequence problem.
\epr

\be{example_matching}
{\rm
Let $H$ be the triangle and let $f_1(x)=f_3(x)=(x-1)^2$ and $f_2(x)=(x-2)^2$.
The graph $L$ constructed in the proof above is given in Figure \ref{figure1}
with nonzero costs of edges indicated in blue. A minimum cost perfect matching
in $L$ is indicated in red, from which we read off an optimal subgraph
$G\subset H$ with edge set $N=\{\{1,2\},\{2,3\}\}$ and optimal value $0$.}
\ee

\section{Unbalanced bipartite graphs}
\label{ubg}

Here we consider our problem over unbalanced bipartite graphs, with one side fixed and small.

\vskip.2cm\noindent{\bf Theorem \ref{unbalanced_bipartite}}
{\it
For any fixed $r$, the optimization problem over degree sequences can be solved in polynomial
time for any given bipartite graph $H\subseteq K_{r,n-r}$ and any given functions $f_1,\dots,f_n$.
}
\vskip.2cm
\bpr
Let $H=(I,J,E)$ and assume that $I=[r]$ and $J=\{r+1,\dots,r+s\}$ with $s:=n-r$.
To find an optimal $G\subseteq H$ we use dynamic programming, namely, we reduce
the problem to computing a shortest path in an auxiliary digraph $D$ with arc lengths,
with nodes labeled by pairs $[v,j]$ with $v\in\{0,1,\dots,s\}^r$ and $j\in\{0,1,\dots,s+1\}$.
More precisely, its node set is
$$\{[0,0]\}\ \uplus\ \{\left[v=(v_1,\dots,v_r),j\right]\ :\
0\leq v_i\leq d_i(H),\ j\in[s]\}\ \uplus\ \{[0,s+1]\}\ .$$
For $j\in[s]$ the node $[v,j]$ (state in the dynamic program) indicates that the current degree of
$i\in[r]=I$, after the neighborhoods of $\{r+1,\dots,r+j\}\subseteq J$ had been decided, is $v_i$.

For $j\in[s]$ there is an arc from $[u,j-1]$ to $[v,j]$ in $D$ (decision in the
dynamic program) if $v-u$ indicates an admissible neighborhood in $G$ for $r+j\in J$,
that is, if $v_i-u_i\in\{0,1\}$ for all $i$ and $v_i-u_i=0$ if $[i,r+j]\notin E$.
Deciding to go along such an arc indicates that the degree of $r+j\in J$ in $G$ will be
$d_{r+j}(G)=\sum_{i=1}^r(v_i-u_i)$, incurring a cost of $f_{r+j}(\sum_{i=1}^r(v_i-u_i))$,
and so we define the length of such an arc to be that cost.
For each $v$, the node (state) $[v,s]$ indicates that after all neighborhoods in $G$
of vertices $r+j\in J$ had been decided, the degree in $G$ of each $i\in I$ is $d_i(G)=v_i$,
and so the arc from $[v,s]$ to $[0,s+1]$ in $D$ is included, with length
set to be the total cost $\sum_{i=1}^r f_i(v_i)$ incurred by all vertices in $I$.

\begin{figure}[hbt]
\centerline{\includegraphics[scale=0.6]{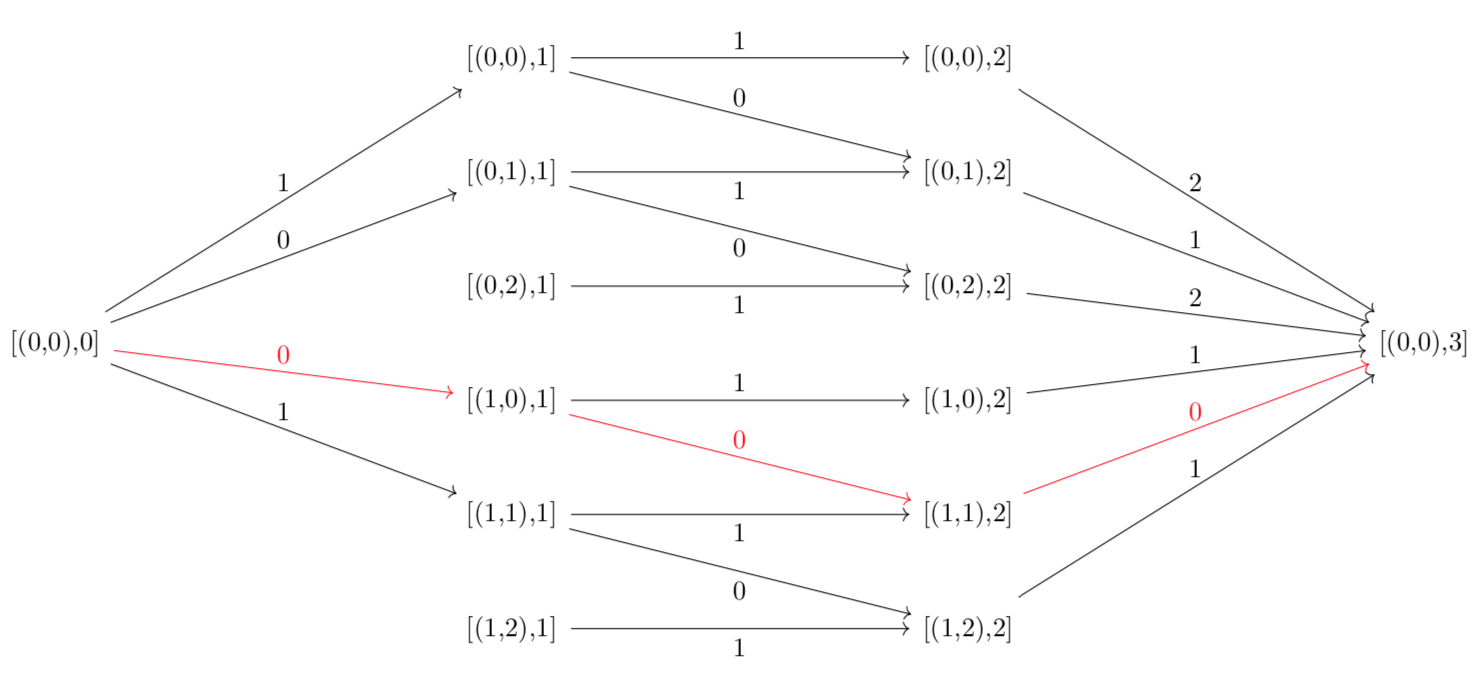}}
\caption{The digraph $D$ of the dynamic program in Example \ref{example}}
\label{figure2}
\end{figure}

Let $[0,0],[v^1,1],\dots,[v^s,s],[0,s+1]$ be a shortest path from $[0,0]$ to $[0,s+1]$ in $D$.
Then it is clear that an optimal subgraph $G=(I,J,F)\subseteq H$
can be read off form it by setting
$$F\ :=\ \{(i,r+j)\ :\ \ i\in I,\ \ r+j\in J,\ \ v^j_i-v^{j-1}_i=1\}\ .$$
So a shortest path in $D$ indeed gives an optimal $G$. Now the number of nodes of $D$ is
$$2+s\prod_{i=1}^r\left(d_i(H)+1\right)\ =\ O\left(n^{r+1}\right)$$
and is polynomial in $n$. Since a shortest path can be found in linear time we are done.
\epr

\be{example}
{\rm
Let $H=(I,J,E)$ be a given bipartite graph with $I=\{1,2\}$, $J=\{3,4\}$, and
$E=\{(1,3),(2,3),(2,4)\}$, and let $f_i(x)=(x-1)^2$ for $i=1,2,3,4$ be given functions.
Figure \ref{figure2} illustrates the digraph $D$ in the construction of
Theorem \ref{unbalanced_bipartite} with edge lengths indicated.
A shortest path in $D$ is indicated in the figure in red, from which we read off an optimal
subgraph $G=(I,J,F)$ with $F=\{(1,3),(2,4)\}$ which is the unique perfect matching in $H$.}
\ee

\section{Monotone functions}
\label{mf}

Here we consider the situation when most functions are all nondecreasing or all nonincreasing.

\vskip.2cm\noindent{\bf Theorem \ref{fixed_monotone}}
{\it
For any fixed $r$, our problem is polynomial time solvable for any graph $H$ and any functions
$f_1,\dots,f_r$ when $f_{r+1},\dots,f_n$ are either all nondecreasing or all nonincreasing.
}
\vskip.2cm
\bpr
Let $H=([n],E)$ and let $I:=[r]$ and $J:=\{r+1,\dots,r+s\}$ with $s:=n-r$. For any
$K\subseteq[n]$ let ${K\choose2}:=\{\{i,j\}:i,j\in K,i\neq j\}$. Suppose first $f_{r+1},\dots,f_n$
are nondecreasing. Let $G=([n],F)\subseteq H$ be an optimal subgraph, that is, attaining minimum
$\sum_{i=1}^n f_i(d_i(G))$, and among such, with $q(G):=|F\cap{J\choose2}|$ minimal.
If there is some $e\in F\cap{J\choose2}$ then removing it from $F$ results in a subgraph $L$
with $\sum_{i=1}^n f_i(d_i(L))\leq\sum_{i=1}^n f_i(d_i(G))$ so $L$ is also optimal with
$q(L)<q(G)$ contradicting the choice of $G$. So $F\cap{J\choose2}=\emptyset$. So our problem
reduces to finding a subgraph $G=([n],F)$ minimizing $\sum_{i=1}^n f_i(d_i(G))$
with $F\cap{J\choose2}=\emptyset$. Now suppose $f_{r+1},\dots,f_n$ are nonincreasing.
Let $G=([n],F)$ be an optimal subgraph and among such, with $q(G):=|F\cap{J\choose2}|$ maximal.
If there is some $e\in E\cap{J\choose2}$ not in $F$ then adding it to $F$ results in a subgraph
$L$ with $\sum_{i=1}^n f_i(d_i(L))\leq\sum_{i=1}^n f_i(d_i(G))$ so $L$ is also optimal with
$q(L)>q(G)$ contradicting the choice of $G$. So $F\cap{J\choose2}=E\cap{J\choose2}$.
Now, for $j\in J$ let $d_j(H[J])$ be its degree in the induced subgraph $H[J]$, and replace the
function $f_j$ by ${\bar f}_j(x):=f_j(x+d_j(H[J]))$ for $x\in\{0,\dots,d_j(H)-d_j(H[J])\}$,
and for $i\in I$ set ${\bar f}_i:=f_i$. Then our problem reduces again to finding a subgraph
$G=([n],F)$ minimizing $\sum_{i=1}^n{\bar f}_i(d_i(G))$ with $F\cap{J\choose2}=\emptyset$.

To find an optimal subgraph $G=([n],F)\subseteq H$ with $F\cap{J\choose2}=\emptyset$,
with respect to the $f_i$ or the ${\bar f_i}$, we use the construction in the proof of
Theorem \ref{unbalanced_bipartite}, applied to the bipartite graph ${\bar H}:=(I,J,{\bar E})$
with ${\bar E}:=\{\{i,j\}\in E: i\in I, j\in J\}$, and extend it so as to allow also decisions
about the existence of edges $\{i,j\}\in E\cap{I\choose2}$, as follows.
Let $t:=|E\cap{I\choose2}|$ and let $\{i_1,j_1\},\dots,\{i_t,j_t\}$ be any
ordering of the edges in $E\cap{I\choose2}$. We add to $D$ the set of nodes
$$\left\{\left[v=(v_1,\dots,v_r),\{i_l,j_l\}\right]\ :\
0\leq v_i\leq d_i(H),\ \ 1\leq l\leq t\right\}\ .$$
We remove the arcs from nodes $[v,s]$ to $[0,s+1]$. We add arcs from
$[v,s]$ to $[v,\{i_1,j_1\}]$ indicating that edge $\{i_1,j_1\}$ is not selected,
and from $[v,s]$ to $[v+{\bf 1}_{i_1}+{\bf 1}_{j_1},\{i_1,j_1\}]$, where ${\bf 1}_l$
is the $l$-th unit vector in $\R^r$, indicating that edge $\{i_1,j_1\}$ is selected.
Also, for $l=2,\dots,t$, we add arcs from $[v,\{i_{l-1},j_{l-1}\}]$ to $[v,\{i_l,j_l\}]$
indicating that edge $\{i_l,j_l\}$ is not selected, and from
$[v,\{i_{l-1},j_{l-1}\}]$ to $[v+{\bf 1}_{i_l}+{\bf 1}_{j_l},\{i_l,j_l\}]$
indicating that edge $\{i_l,j_l\}$ is selected. To all these arcs we assign length zero.
Finally, for each $v$, we add an arc from $[v,\{i_t,j_t\}]$ to $[0,s+1]$, with length
set to be the total cost $\sum_{i=1}^r f_i(v_i)$ incurred by all vertices in $I$. Let
$$[0,0],[v^1,1],\dots,[v^s,s],
[v^{i_1,j_1},\{i_1,j_1\}],\dots,[v^{i_t,j_t},\{i_t,j_t\}],[0,s+1]$$
be a shortest path from $[0,0]$ to $[0,s+1]$ in $D$.
Then an optimal graph is $G=([n],F)$ with
\begin{eqnarray*}
F & :=& \{\{i,r+j\}\ :\ \ i\in I,\ \ r+j\in J,\ \ v^j_i-v^{j-1}_i=1\} \\
& \uplus & \ \{\{i_1,j_1\}\ :\ v^{i_1,j_1}-v^s={\bf 1}_{i_1}+{\bf 1}_{j_1}\}\\
& \uplus & \ \{\{i_l,j_l\}\ :\
v^{i_l,j_l}-v^{i_{l-1},j_{l-1}}={\bf 1}_{i_l}+{\bf 1}_{j_l},\ \ 2\leq l\leq t \}.
\end{eqnarray*}

Now the number of nodes of $D$ is
$$2+\left(s+t\right)\prod_{i=1}^r\left(d_i(H)+1\right)\ =\ O\left(n^{r+1}\right)$$
and is polynomial in $n$. Since a shortest path can be found in linear time we are done.
\epr

\section*{Acknowledgments}
G. Deza was partially supported by an internship from the university of Toronto. S. Onn was
partially supported by a grant from the Israel Science Foundation and by the Dresner chair.

\section*{Appendix}

\begin{figure}[hbt]
\centerline{\includegraphics[scale=0.6]{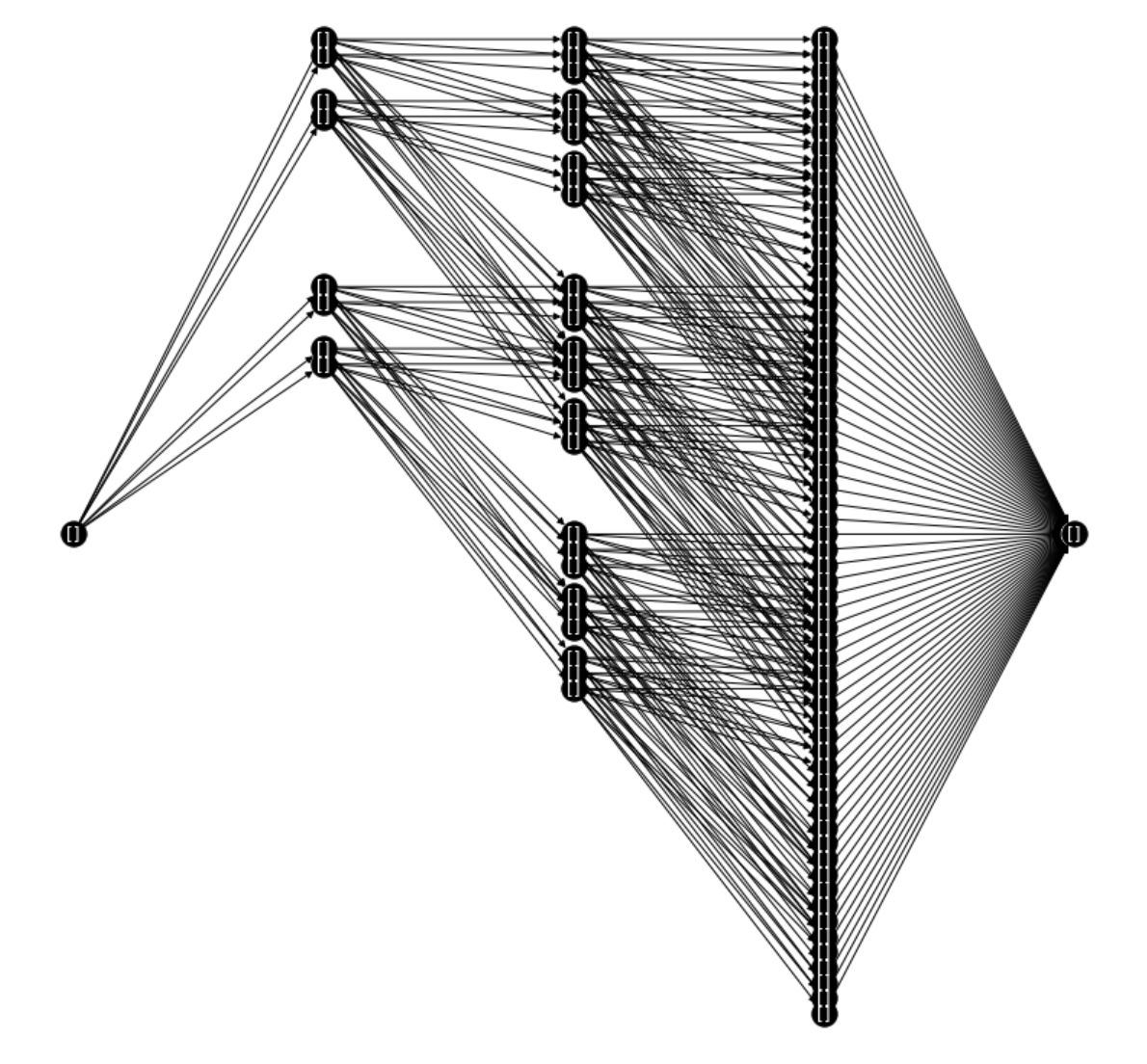}}
\caption{The digraph $D$ of the dynamic program applied to $H:=K_{3,3}$}
\label{figure3}
\end{figure}

The algorithm of Theorem \ref{unbalanced_bipartite} has been implemented and
is available in \cite{Dez}. For instance, when applied to $H:=K_{3,3}$, the
directed graph $D$ of the dynamic program that the algorithm produces is depicted
in Figure \ref{figure3} (with all vertices unreachable from the source removed).

\end{document}